\documentclass{amsproc}
\usepackage[matrix,arrow,curve]{xy}
\usepackage{amssymb,amsthm}
\usepackage{supertabular}

\newcommand{\ddate}{February 16, 2007}

\newtheorem{dummy}{anything}[section]

\newtheorem{Lemma}[dummy]{Lemma}
\newtheorem{Proposition}[dummy]{Proposition}

\newtheorem{Example}[dummy]{Example}

\newtheorem{Remark}[dummy]{Remark}

\newtheorem{ccote}[dummy]{}

\newcommand{\bbr}{{\mathbb R}}

\newcommand{\bbc}{{\mathbb C}}
\newcommand{\bbz}{{\mathbb Z}}

\newcommand{\calc}{{\mathcal C}}
\newcommand{\cald}{{\mathcal D}}

\newcommand{\calh}{{\mathcal H}}

\newcommand{\calk}{{\mathcal K}}

\newcommand{\caln}{{\mathcal N}}

\newcommand{\cals}{{\mathcal S}}

\newcommand{\llangle}[1]{\langle #1 \rangle}
\newcommand{\pcirc}{\kern .7pt {\scriptstyle \circ} \kern 1pt}

\newcommand{\mun}{{-1}}

\newcommand{\rmpos}{\ensuremath{(\bbr_{>0})^m}}
\newcommand{\rmcr}{\ensuremath{\bbr_{\scriptscriptstyle\nearrow}^m}}

\newcommand{\fpp}{\ensuremath{\hookrightarrow}}
\newcommand{\chrmcr}{{\rm Ch}(\bbr_{\scriptscriptstyle\nearrow}^m)}
\newcommand{\chrmpos}{{\rm Ch}( (\bbr_{>0})^m)}
\newcommand{\ch}{{\rm Ch}}
\newcommand{\nua}[2]{\caln^{#1}_{#2}}

\newcommand{\cha}[2]{\calc^{#1}_{#2}}

\newcommand{\psetminus}{{\scriptstyle\setminus}}

\renewcommand{\:}{\colon}

\newcommand{\ccvi}[4]{  & {#4} & {#3} & {#2}\\[1.3mm]}
\newcommand{\pvi}[7]{\stepcounter{compol}{\small\thecompol} & %
$\langle #7\rangle$ & $(#4)$ }
\newcommand{\fno}[1]{\, $(^{#1})$}

\newcommand{\preu}{\noindent{\sc Proof: \ }}
\newcommand{\sk}[1]{\vskip #1 mm}
\newcommand{\equref}[1]{(\ref{#1})}
\newcommand{\hfl}[2]{\smash{\mathop{\hbox to 1 truecm{\kern %
3pt\rightarrowfill\kern 3pt}}%
\limits^{\scriptstyle#1}_{\scriptstyle#2}}}
\newcommand{\cqfd}{\unskip\kern 6pt\penalty 500%
\raise -2pt\hbox{\vrule\vbox to10pt{\hrule width %
4pt\vfill\hrule}\vrule}\smallskip}
\newcommand{\proref}[1]{Proposition~\ref{#1}}
\newcommand{\remref}[1]{Remark~\ref{#1}}
\newcommand{\lemref}[1]{Lemma~\ref{#1}}

\newcommand{\exref}[1]{Example~\ref{#1}}

\newcounter{compol}

\newcommand{\undm}[1]{\{1,\dots,#1\}}

\title{Geometric descriptions of polygon and chain spaces}
\author{Jean-Claude HAUSMANN}
\date{\ddate}
\address{Section de Math\'ematiques, Universit\'e de Gen\`eve,
\\ B.P. 240, CH-1211 Geneva 24, Switzerland}
\email{hausmann@math.unige.ch}
\subjclass{Primary 55R80, 70G40; Secondary 57R65}

\begin{document}
\maketitle 

\begin{abstract}
We give a few simple methods to geometically describe some polygon and chain spaces
in $\bbr^d$. They are strong enough to give tables of $m$-gons and $m$-chains
when $m\leq 6$.
\end{abstract}

\section*{Introduction}

For $a=(a_1,\dots,a_m)\in\bbr_{>0}^m$ and $d$ an integer, define the subspace
$\cha{m}{d}(a)$ of $\prod_{i=1}^{m-1}S^{d-1}$ by
$$
\cha{m}{d}(a)=\big\{z=(z_1,\dots,z_{m-1})\in\prod_{i=1}^{m-1}S^{d-1}
\mid \sum_{i=1}^{m-1}a_iz_i=a_m\, e_1
\big\} \, ,
$$
where $e_1=(1,0,\dots,0)$ is the first vector of the standard
basis $e_1,\dots,e_d$ of $\bbr^d$.
An element of $\cha{m}{d}(a)$, called a {\it chain}, can be visualized
as a configuration of $(m-1)$-segments in $\bbr^d$, of length
$a_1,\dots,a_{m-1}$, joining the origin to $a_me_1$.
The group $O(d-1)$,
seen as the subgroup of $O(d)$ stabilizing the first axis, acts naturally
(on the left) on
$\cha{m}{d}(a)$. The quotient space by $SO(d-1)$ coincides with the
{\it polygon space}
$$
\begin{array}{rcl}
\nua{m}{d}(a) &=&
SO(d-1)\big\backslash \cha{m}{d}(a)
\\ &\approx&
SO(d)\bigg\backslash\big\{
\rho=(\rho_1,\dots\rho_m)\in (\bbr^d)^m \mid
|\rho_i|=a_i \hbox{ and } \sum_{i=1}^m\rho_i=0
\big\} \ .
\end{array}
$$
The notations are that of \cite{HR} where it is emphasized how the union  $\nua{m}{d}$
of $\nua{m}{d}(a)$ for all $a\in\bbr_{>0}^m$ is related to the spaces studied in
statistical shape analysis (see, e.g.~\cite{KBCL}).
An element $a\in\bbr_{>0}^m$ is {\it generic}
if $\cha{m}{1}(a)=\emptyset$, that is to say there is no lined chain or polygon
configuration.
When $a$ is generic, $\cha{m}{d}(a)$ is a smooth closed
manifold of dimension $(m-2)(d-1)-1$ (see, e.g.~\cite{Ha}).

Mathematical robotics is specially interested in the chain
and polygon spaces for $d=2,3$.
When $a$ is generic, the action of $SO(d-1)$ on $\cha{m}{d}$ is then free
and therefore $\nua{m}{2}(a)$ and $\nua{m}{3}(a)$ are closed smooth manifolds of
dimension $m-3$ and $2(m-3)$ respectively (in addition, $\nua{m}{3}(a)$
carries a symplectic structure, see e.g. \cite{KM2}).
One has $\cha{m}{2}(a)=\nua{m}{2}(a)$ and $\cha{m}{3}(a)\to\nua{m}{3}(a)$ is a
principal circle bundle.

In this paper, we present a few geometrical methods permitting us to describe in some cases the spaces $\cha{m}{d}(a)$ and $\nua{m}{d}(a)$. From the classification
results (see Section~\ref{Sclassi}),
this enables us to describe all the chain or polygon spaces in $\bbr^d$ when $m\leq 6$
(tables in Section~\ref{Stables}).

\section{Review of the classification results}\label{Sclassi}

The idea of the classification of the polygon and chain spaces goes back to \cite{Wa}.
Details may be found in \cite{HR}.

\begin{ccote} Short subsets. \rm
Let $a=(a_1,\dots,a_m)\in\bbr_{>0}^{m}$. A subset $J$ of
$\undm{m}$ is called {\it short} if $\sum_{i\in J}a_i<\sum_{i\notin J}a_i$.
Short subsets form, with inclusion, a poset $\cals(a)$. Define
$\cals_m(a)=\{J\in\cals(a)\mid m\in J\}$.

\begin{Lemma}\label{shortdetcha}
Let $a$ and $a'$ be generic elements in $\bbr_{>0}^m$.
Suppose that $\cals_m(a)$ and $\cals_m(a')$ are poset isomorphic.
Then:
\begin{enumerate}
\renewcommand{\labelenumi}{(\roman{enumi})}
\item $\cha{m}{d}(a)$ and $\cha{m}{d}(a')$ are $O(d-1)$-equivariantly
diffeomorphic.
\item $\nua{m}{d}(a)$ and $\nua{m}{d}(a')$ are diffeomorphic.
\end{enumerate}
\end{Lemma}

\preu
If $\cals_m(a)\approx\cals_m(a')$, then
there is a poset isomorphism $\varphi\:\cals(a)\stackrel{\approx}{\to}\cals(a')$
with $\varphi(m)=m$
(see \cite[Proposition~2.5]{HK}).
It is well known that $\cals(a)\approx\cals(a')$
implies (ii) (see, e.g. \cite[Proposition~2.2]{HK} or
\cite[Theorem~1.1]{HR}). We give however the variation of the proof to get the
less classical (stronger) fact that $\cals(a)\approx\cals(a')$
implies (i).

Let $\calk_d(a)=\{z=(z_1,\dots,z_m)\in\prod_{i=1}^m S^{d-1}\mid \sum_{i=1}^ma_iz_i=0\}$.
The group $O(d)$ acts on the left on $\calk_d(a)$ and
$\nua{m}{d}(a)=SO(d)\backslash\calk_d(a)$.
The function $F:\calk_d(a)\to S^{d-1}$ given by
$F(z)=z_m$ is a submersion
(since $F$ is $O(d)$-equivariant). One has $\cha{m}{d}(a)=F^\mun(-e_1)$ with
its residual $O(d-1)$-action.

Let $\sigma$ be the bijection of $\undm{m-1}$ giving the poset
isomorphism $\cals_m(a)\stackrel{\approx}{\to}\cals_m(a')$ and then
$\cals(a)\stackrel{\approx}{\to}\cals(a')$.
Then
$(z_1,\dots,z_{m-1},z_m)\mapsto (z_{\sigma(1)},\dots,z_{\sigma(m-1)},z_m)$ induces
a $O(d-1)$-equivariant diffeomorphism from $\cha{m}{d}(a_1,\dots,a_{m-1},a_m)$ onto 
$\cha{m}{d}(a_{\sigma(1)},\dots,a_{\sigma(m-1)},a_m)$. We can therefore suppose
that $\cals(a)=\cals(a')$ and $\sigma={\rm id}$. We claim that
$\cha{m}{d}(a)$ and $\cha{m}{d}(a')$ are then canonically diffeomorphic. 
Indeed, if $\cals(a)=\cals(a')$, the segment $[a,a']$ contains only generic elements.
Hence, the union 
$$
X=\bigcup_{b\in [a,a']}\big(\calk_d(b)\times\{b\}\big) \subset 
\big(\prod_{i=1}^m S^{d-1}\big)\times [a,a']
$$ 
is an $O(d)$-cobordism
between $\calk_d(a)$ and $\calk(a')$ and the projection $\pi:X\to [a,a']$
has no critical point. One still has the map $F:X\to S^{d-1}$
given by $F(z,t)=z_m$ and $Y=F^\mun(-e_1)$ is an
$O(d-1)$-cobordism between $\cha{m}{d}(a)$ and $\cha{m}{d}(a')$, with
again the projection $\pi$ over $[a,a']$ being a submersion. The standard metric 
on $\prod_{i=1}^{m-1}S^{d-1}$ induces an $O(d-1)$ invariant Riemannian metric on $Y$. 
Following the gradient lines of $\pi$ for this metric gives the required
$O(d-1)$-equivariant diffeomorphism 
$\Psi:\cha{m}{d}(a)\stackrel{\approx}{\to}\cha{m}{d}(a')$.
\cqfd
\end{ccote}

\begin{ccote} Walls and chambers.  \rm
For $J\subset\undm{m})$,
let $\calh_J$ be the hyperplane ({\it wall}) of $\bbr^m$ defined by
$$
\calh_J:=\Big\{(a_1,\dots,a_m)\in\bbr^m \Bigm|
\sum_{i\in J}a_i=\sum_{i\notin J}a_i\Big\}.
$$
The union $\calh(\bbr^m)$ of all these walls
determines a set $\chrmpos$ of open {\it chambers} in $\rmpos$
whose union is the set of generic elements.
Two generic elements $a$ and $a'$ are in the same chamber if and only if
$\cals(a)=\cals(a')$. We call $\ch(a)$ the chamber of a generic element $a$.
If $\alpha$ is a chamber, the poset $\cals(a)$ is the same for all $a\in\alpha$
and is denoted by $\cals(\alpha)$.
\end{ccote}

\begin{ccote} Permutations.  \rm
Let $\sigma$ be a permutation of $\undm{m}$. The map which sends
$(z_1,\dots,z_m)$ to $(z_{\sigma(1)},\dots,z_{\sigma(m)})$ induces
a diffeomorphism from $\nua{m}{d}(a_1,\dots,a_m)$ onto
$\nua{m}{d}(a_{\sigma(1)},\dots,a_{\sigma(m)})$. For the sake of the
classification of $\nua{m}{d}(a)$, we may as well assume that $a\in\rmcr$
where
$$\rmcr:=\{(a_1,\dots,a_m)\in\bbr^m\mid
0< a_1\leq\cdots\leq a_m\} \, .
$$
Observe that we then do not classify all the chain spaces $\cha{m}{d}(a)$
but only those for which $a_m\geq a_i$ for $i<m$. Indeed,
the permutation $\sigma$ induces a diffeomorphism from
$\cha{m}{d}(a_1,\dots,a_m)$ onto
$\cha{m}{d}(a_{\sigma(1)},\dots,a_{\sigma(m)})$ if and only if $\sigma(m)=m$.
We denote by $\chrmcr$ the set of chambers determined in $\rmcr$ by 
the hyperplane arrangement $\calh(\bbr^m)$.
\end{ccote}

\begin{ccote} The genetic code of a chamber.  \rm
A chamber $\alpha\in\chrmcr$ is determined by $\cals(\alpha)$ which, in turn,
is determined by $\cals_m(\alpha)$. 
Consider the partial order ``\fpp'' on
the subsets of $\undm{m}$ where $A\fpp B$ if and only if there exits
a non-decreasing map $\varphi:A\to B$ such that $\varphi(x)\geq x$.
For instance $X\fpp Y$ if $X\subset Y$ since one can
take $\varphi$ being the inclusion.
The {\it genetic code} of $\alpha$ is the set of
elements $A_1,\dots,A_k$ of $S_m(\alpha)$ which are
maximal with respect to the order ``\fpp''.
Thus, the chamber $\alpha$ is determined
by its genetic code; we write
$\alpha=\llangle{A_1,\dots,A_k}$ and call the sets $A_i$
the {\it genes} of $\alpha$. As, in this paper $m\leq 9$, we abbreviate
a subset $A$ by the sequence of its digits, e.g. $\{6,2,1\}=621$.
In \cite{HR}, an algorithm is presented to list by their genetic codes
all the elements of  $\chrmcr$ and then all the chambers up to permutation of the components.
Tables for $m\leq 6$ are given in \cite{HR} (and in Section~\ref{Stables} below); more tables, for $m\leq 9$, may be found in \cite{HRWeb}. 
The algorithm produces, in each chamber $\alpha$,
a representative $a_{\rm min}(\alpha)\in\alpha$;
though this is not proved theoretically, $a_{\rm min}(\alpha)$
turned out in all known cases to have 
integral components $a_i$ and minimal
$\sum a_i$. See examples in the tables below.
\end{ccote}

\section{Procedures of description}

\subsection{Adding a tiny edge}\label{atiny}

Let $a=(a_2,\dots,a_m)$ be a generic element of $\bbr_{\scriptscriptstyle\nearrow}^{m-1}$.
If $\varepsilon>0$ is small enough, the $m$-tuple
$a^+:=(\delta,a_2,\dots,a_m)$ is a generic
element of $\rmcr$ for $0<\delta\leq\varepsilon$. This defines a
map $\ch(\bbr_{\scriptscriptstyle\nearrow}^{m-1})\stackrel{+}{\to}\chrmcr$,
sending $\alpha$ to $\alpha^+$, which is injective (see \cite[Lemma~5.1]{HR}).
The genetic code of $\alpha^+$ has the same number of genes
than that of $\alpha$ and the correspondence goes as follows.
If $\{p_1,\dots, p_r\}$ is a gene of $\alpha$,
then $\{p_1^+,\dots, p_r^+,1\}$ is a gene of $\alpha^+$,
where $p_i^+=p_i+1$. For example: $\llangle{631,65}^+=\llangle{7421,761}$.
The minimal integral representative $a_{\rm min}(\alpha^+)$ of $\alpha^+$ is a {\it conventional representative}: it starts with a $0$ followed by the components
of $a_{\rm min}(\alpha)$. Example: as $a_{\rm min}(\llangle{3})=(1,1,1)$,
then $a_{\rm min}(\llangle{3}^+)=a_{\rm min}(\llangle{41})=(0,1,1,1)$,
$a_{\rm min}(\llangle{41}^+)=a_{\rm min}(\llangle{521})=(0,0,1,1,1)$, {\it etc}.
It has to be understood that these
vanishing components stand for small enough positive
real numbers, whose sum is less than $1$.   

\begin{Proposition}\label{Paddtiny}
There is a $O(d-1)$-equivariant diffeomorphism
$$
\Phi\:\cha{m}{d}(\alpha^+)\stackrel{\approx}{\longrightarrow} S^{d-1}\times\cha{m-1}{d}(\alpha) \, ,
$$
where $S^{d-1}\times\cha{m}{d-1}(\alpha)$ is equipped with the diagonal $O(d-1)$-action.
\end{Proposition}

\preu
Let $a=(a_2,\dots,a_m)\in\alpha$ and $a^+=(\varepsilon,a_2,\dots,a_m)\in\alpha^+$.
The map $\Phi$ is of the form $(\Phi_1,\Phi_2)$, where
$\Phi_1\:\cha{m}{d}(a^+)\to S^{d-1}$ and
$\Phi_2\:\cha{m}{d}(a^+)\to \cha{m-1}{d}(a)$
are $O(d-1)$-equivariant maps.
The map $\Phi_1$ is just given by $\Phi_1(z_1,\dots,z_m)=z_1$.
It remains to define $\Phi_2$.

If $p\in\bbr^d$ satisfies $p\neq -|p|e_1$, there is a unique $R_p\in SO(d)$ such that
$R_p(p)=|p|e_1$ and $R_p(q)=q$ if $q\in{\rm EV}(p,e_1)^\perp$, the orthogonal complement to the
vector space ${\rm EV}(p,e_1)$ generated by $p$ and $e_1$.
In particular, $R_{e_1}={\rm id}$. The map $p\to R_p$ is smooth. We shall apply that to
$p=p(z)$, where
$$
p(z)=\sum_{i=2}^m a_iz_i = a_m e_1 - \varepsilon z_1 \, .
$$
We may suppose that $\varepsilon<a_m$, so $p(z)\neq -|p(z)|e_1$.
The correspondence $(z_1,\dots,z_m)\mapsto (R_{p(z)}z_2,\dots,R_{p(z)}z_m)$ gives a smooth map
$$
\Phi_2'\:\cha{m}{d}(a^+)\to\cha{m-1}{d}(a_2,\dots,a_{m-1},|p(z)|) \, .
$$
 The fact that
$(\delta,a_2,\dots,a_m)$ is generic when $0<\delta\leq\varepsilon$ implies that
$$
\ch(a_2,\dots,a_{m-1},|p(z)|)=\ch(a) \, .
$$
We can then use the canonical $O(d-1)$-equivariant diffeomorphism
$$
\Psi\:\cha{m-1}{d}(a_2,\dots,a_{m-1},|p(z)|)\stackrel{\approx}{\to}
\cha{m-1}{d}(a_2,\dots,a_m)
$$ 
constructed in the proof of \lemref{shortdetcha} 
and define $\Phi_2=\Psi\pcirc\Phi_2'$.
If $A\in O(d-1)$, the formula $A R_{p}= R_{A(p)}A$ holds in $O(d-1)$, as easily seen
on ${\rm EV}(p,e_1)$ and on ${\rm EV}(p,e_1)^\perp$. This implies that
$\phi_2'$ is $O(d-1)$-equivariant.

We have thus constructed an $O(d-1)$-equivariant smooth map
$\Phi\:\cha{m}{d}(\alpha^+)\to S^{d-1}\times\cha{m-1}{d}(\alpha)$.
The reader will easily figure out what the inverse $\Phi^\mun$ of $\Phi$ is like, proving that
$\Phi$ is a diffeomorphism. \cqfd

We now turn our interest to $\nua{m}{3}(\alpha^+)$. Let $\cald(\alpha)$ be the total space
of the $D^2$-disk bundle associated to $\cha{m-1}{3}(\alpha)\to\nua{m-1}{3}(\alpha)$.
We call {\it double of } $\cald(\alpha)$
the union of two copies of $\cald(\alpha)$,
with opposite orientations, along their common boundary $\cha{m}{3}(\alpha)$.

\begin{Proposition}\label{Paddtiny-nua}
\renewcommand{\labelenumi}{(\alph{enumi})}
\hskip -10mm
\begin{enumerate}
\item  $\nua{m}{3}(\alpha^+)$ is diffeomorphic to 
$S^2\times_{S^1}\cha{m-1}{3}(\alpha)$.
\item  $\nua{m}{3}(\alpha^+)$ is diffeomorphic to the double of $\cald(\alpha)$.
\end{enumerate}
\end{Proposition}

In Part (a), $S^2\times_{S^1}\cha{m-1}{3}(\alpha)$ denotes the quotient of $S^2\times\cha{m-1}{3}(\alpha)$
by the diagonal action of $S^1=SO(2)$. The projection $S^2\times_{S^1}\cha{m-1}{3}(\alpha)\to\nua{m-1}{3}(\alpha)$
is then the $S^2$-associated bundle to the $SO(2)$-principal bundle $\cha{m-1}{3}(\alpha)\to\nua{m-1}{3}(\alpha)$.
A direct proof of Part (b) may be found in \cite[Prop.\,6.4]{HR}.

\preu
For Part (a), we check that the diffeomorphism
$$
\Phi\:\cha{m}{3}(\alpha^+)\stackrel{\approx}{\to} S^2\times\cha{m-1}{3}(\alpha) 
$$
of \proref{Paddtiny} descends to a diffeomorphism from $\nua{m}{3}(\alpha^+)$ to
$S^2\times_{S^1}\cha{m-1}{3}(\alpha)$. For Part(b), we observe that
$\cald(\alpha)$ is the mapping cylinder of the
projection $\cha{m-1}{3}(\alpha)\to\nua{m-1}{3}(\alpha)$. The double of
$\cald(\alpha)$ is then diffeomorphic to
$M=[-1,1]\times \cha{m-1}{3}(\alpha)\big/\!\sim$, where ``$\sim$''
is the equivalence relation generated by $(-1,z)\sim (-1,Az)$ and
$(1,z)\sim (1,Az)$ for all $z\in\cha{m-1}{3}(\alpha)$ and all $A\in S^1$.
Each $S^1$-orbit of $S^2$ has an unique point of the form  $(u_1,0,u_3)$.
To $(u,z)\in S^2\times\cha{m-1}{d}(\alpha)$ with $u=(u_1,0,u_3)$, we associate
the class $[u_1,z]$ in $M$ and check that this
correspondence gives rise to a diffeomorphism from
$S^2\times_{S^1}\cha{m-1}{3}(\alpha)$ to the double of $\cald(\alpha)$.
\cqfd

\begin{Example}\label{111}\rm
When $m=3$, there is only one chamber $\alpha=\llangle{3}$,
with $a_{\rm min}(\alpha)=(1,1,1)$, for which $\cha{3}{d}(\alpha)$ is not empty.
Its image under adding tiny edges gives a chamber
$\llangle{\{m,m-3,m-2,\dots,1\}}\in\ch(\bbr_{\scriptscriptstyle\nearrow}^m)$ with
$a_{\rm min}=(0,\dots,0,1,1,1)$ (conventional representative, \S~\ref{atiny}). As $\cha{3}{d}\llangle{3}=S^{d-2}$
with the standard $O(d-1)$-action,
Propositions~\ref{Paddtiny} and~\ref{Paddtiny-nua} give the following

\begin{center}
\begin{tabular}{c|c|c|c}
\multicolumn{4}{c}{\bf The chamber $\alpha=\llangle{\{m,m-3,m-2,\dots,1\}}$ }
\\ \hline\rule{0pt}{3ex}
$a_{\rm min}(\alpha)$ & $\nua{m}{2}(\alpha)$ & $\nua{m}{3}(\alpha)$ & $\cha{m}{d}(\alpha)$
\\[1pt]
\hline\rule{0pt}{3ex}
$(0,\dots,0,1,1,1)$ &  $T^{m-3}{\scriptstyle\coprod}T^{m-3}$  &  $(S^2)^{m-3}$ &
$(S^{d-1})^{m-3}\times S^{d-2}$
\\ \hline
\end{tabular}
\end{center}
\end{Example}

\begin{Remark}\label{111-rem}\rm
Let $A=\{m,m-3,m-2,\dots,1\}$. We claim that
$\alpha=\llangle{A}$ as above is the only chamber in $\rmcr$ having
$J\in\cals_m(\alpha)$ with $|J|=m-3$. Indeed, let $\beta\in\chrmcr$ having
$J\in\cals_m(\beta)$ with $J\neq A$ and $|J|=m-3$.
Then $A'=\{m,m-2,m-4,\dots,1\}$ would satisfy $A'\fpp J$.
Then, $\bar A'=\{m-3,m-1\}$ would be long, which contradicts
$\{m-3,m-1\}\fpp\{m-2,m\}\in\cals_m(\beta)$. Now, if $A\in\cals_m(\beta)$,
then $\{m,m-2\}$ is long, since $\overline{\{m,m-2\}}\fpp A$.
Therefore, $A\in\cals_m(\beta)$ implies $\beta=\llangle{A}$. For an application
of this remark, see Propositions~\ref{connectivity} and~\ref{Psurgery2}.
\end{Remark}

\subsection{The manifold $V_d(a)$}

Let $a\in\bbr_{>0}^m$. Define
$$
V_d(a)=\{z=(z_1,\dots,z_{m-1})\in \prod_{i=1}^{m-1} S^{d-1} \mid
\sum_{i=1}^{m-1}a_iz_i=t e_1 \hbox{ with } t\geq a_m\} \, .
$$
Let $f:V_d(a)\to\bbr$ defined by $f(z)=-|\sum_{i=1}^{m-1}a_iz_i|$.
The group $O(d-1)$ acts on $V_d(a)$.
The following proposition is proven in \cite[Th.~3.2]{Ha}.

\begin{Proposition}\label{Morse}
Suppose that $a\in\bbr_{>0}^m$ is generic. Then
\renewcommand{\labelenumi}{(\roman{enumi})}
\begin{enumerate}
\item  $V_d(a)$ is a smooth $O(d-1)$-submanifold
of $\prod_{i=1}^{m-1} S^{d-1}$, of dimension $(m-2)(d-1)$, with boundary
$\cha{m}{d}(a)$.
\item $f$ is a $O(d-1)$-equivariant Morse function, with one critical point $p_J$
for each $J\in\cals_m(a)$, where
$p_J=(z_1,\dots,z_{m-1})$ with $z_i$ equal to $-e_1$ if $i\in J$ and $e_1$
otherwise (aligned configuration). The index of $p_J$ is $(d-1)(|J|-1)$. \cqfd
\end{enumerate}
\end{Proposition}

This permits us to get some information on $\cha{m}{d}(a)$. 

\begin{Example}\label{mmm} The chamber $\llangle{m}$. \rm
If $\cals_m=\{m\}$, $f\:V_d(a)\to\bbr$ has only one critical point, of index $0$.
Hence, $\cha{m}{d}(a)\approx S^{(m-2)(d-1)-1}$ and the $O(d-1)$ action
is conjugate to that obtained by the embedding 
$S^{(m-2)(d-1)-1}\subset (\bbr^d)^{m-2}$ with the standard 
diagonal action \cite[Prop.\,4.2]{Ha}.  
The chamber of $a$ has here genetic code $\llangle{m}$, with minimal representative
$(1,\dots,1,m-2)$. One then has:

\begin{center}
\begin{tabular}{c|c|c|c}
\multicolumn{4}{c}{\bf The chamber $\alpha=\llangle{m}$ }
\\ \hline\rule{0pt}{3ex}
$a_{\rm min}(\alpha)$ & $\nua{m}{2}(\alpha)$ & $\nua{m}{3}(\alpha)$ & $\cha{m}{d}(\alpha)$
\\[1pt]
\hline\rule{0pt}{3ex}
$(1,\dots,1,m-2)$ &  $S^{m-3}$  &  $\bbc P^{m-3}$ &
$S^{(m-2)(d-1)-1}$
\\ \hline
\end{tabular}
\end{center}
\end{Example}

\sk{1}
Another consequence of \proref{Morse} is the connectivity of $\cha{m}{d}(\alpha)$.
We saw in \exref{111} that $\cha{m}{d}(\beta)=(S^{d-1})^{m-3}\times S^{d-2}$
if $\beta=\llangle{\{m,m-3,m-2,\dots,1\}}$. Thus, $\pi_{d-2}(\cha{m}{d}(\beta))\approx\bbz$ if $d\geq 3$ and 
$\pi_{0}(\cha{m}{2}(\beta))$ has $2$ elements.   
But this is an exceptional case:

\begin{Proposition}\label{connectivity}
Let $\alpha$ be a chamber of $\rmcr$ with $\alpha\neq\llangle{\{m,m-3,m-2,\dots,1\}}$.
Then, $\cha{m}{d}(\alpha)$ is $(d-2)$-connected, i.e.
$\pi_i(\cha{m}{d}(\alpha))=0$ for $i\leq d-2$.
\end{Proposition}

\preu Let $a\in\rmcr$ be a representative of $\alpha$.
By \proref{Morse}, one has that $\pi_i(V_d(a))=0$ if $i\leq d-2$.
If $\alpha\neq\llangle{\{m,m-3,m-2,\dots,1\}}$, then $|J|\leq m-3$
for all $J\in\cals_m(a)$ by \remref{111-rem}.
Then $V_d(a)$ has a handle decomposition,
starting from $\cha{m}{d}(a)$, with handles of index
$\geq (m-2)(d-1)-(m-4)(d-1)=2(d-1)$. Therefore,
$\pi_i(\cha{m}{d}(a))\approx\pi_i(V_d(a))$ for $i\leq 2d-2>d-2$. \cqfd

\begin{Remark}\label{connectivity-rem}\rm
 When $d=2$, \proref{connectivity}
says that $\llangle{\{m,m-3,m-2,\dots,1\}}$ is the only chamber $\beta$ of $\rmcr$
for which $\nua{m}{2}(\beta)$ is not connected. This was proved by
Kapovich and Millson \cite{KM} (see also \cite[Ex.2 in \S 1]{FS}).
\end{Remark}

\subsection{Crossing walls and surgeries}\label{SSurg}

\begin{Proposition}\label{Psurgery1}
Let $J\subset\{1,\dots,m\}$, defining the wall $\calh_J$ in $\bbr^m$.
Let $\alpha$ and $\beta$ be two chambers of $\bbr^m$, with
$\cals_m(\beta)=\cals_m(\alpha)\cup \{J\}$.
Then $\cha{m}{d}(\beta)$
is obtained from $\cha{m}{d}(\alpha)$ by an $O(d-1)$-equivariant
surgery of index $A=(d-1)(|J|-1)-1$:
$$
\cha{m}{d}(\beta)\approx\big(\cha{m}{d}(\alpha)\psetminus\, (S^{A}\times D^{B})\big)
\cup_{S^{A}\times S^{B-1}}
\big(D^{A+1}\times S^{B-1}\big) \, ,
$$
with $B= (m-1-|J|)(d-1)$. The $O(d-1)$-action on $D^{A+1}$ and $D^B$ comes from their natural
embedding into a product of copies of $\bbr^{d-1}$ with the diagonal action.
\end{Proposition}

\proof
Let $a\in\alpha$ and $b\in\beta$. As $\cals_m(\beta)=\cals_m(\alpha)\cup \{J\}$,
the segment $[a,b]$ in $\bbr^m$ crosses the wall $\calh_J$ and has no intersection
with any other wall. There exists a vector orthogonal to $\calh_J$
with coordinates equal to $\pm 1$.
Therefore, $e_1$ is transverse to $\calh_J$ and, by changing $a$ and $b$
if necessary, we assume that $a=b+\lambda e_1$. By \proref{Morse},
the manifold $V_d(b)\psetminus {\rm int} V_d(a)$ is a $O(d-1)$-equivariant
cobordism $W$ from $\cha{m}{d}(a)$ to $\cha{m}{d}(b)$.
The map $f:W\to \bbr$ defined by $f(\rho)=-|\sum_{j=1}^{m-1}a_j\rho_j|$
is an invariant Morse function having a single critical
point $\rho^0$ of index $(d-1)(|J|-1)$; the components $(\rho_1^0,\dots,\rho_{m-1}^0)$ of $\rho^0$ satisfy
$\rho_i^0=-e_1$ if $i\in J$ and $\rho_i^0=e_1$ if $i\notin J$. By relabeling the $\rho_i$ if necessary,
we assume that $J=\{1,2,\dots,k,m\}$ and $\bar J = \{k+1,\dots,m-1\}$. The index
of $\rho^0$ is then equal to $(d-1)k$. Therefore, $W$ is obtained by adding to a collar neighborhood of $\cha{m}{d}(a)$ an $O(d-1)$-equivariant handle of index $(d-1)k$,
whence the surgery assertion.
For a reference about equivariant Morse theory, see \cite[\S~4]{Wn}.
By \cite[Lemma~4.5]{Wn}, the $O(d-1)$-action is determined by the linear isotropy action on $T_{\rho^0}W$, which we shall now describe.

Let $K_\rho=\sum_{j=1}^{k}a_j\rho_j$ and $L_\rho=K_\rho-a_me_1$.
Let $p_1:\bbr^d\to\bbr^{1}$ and $P:\bbr^d\to\bbr^{d-1}$
be the maps $p_1(x_1,\dots,x_d)=x_1$ and $P(x_1,\dots,x_d)=(x_2,\dots,x_d)$
For $\varepsilon >0$, we consider the following open neighborhood $\caln_\varepsilon$ of $\rho^0$ in $W$
$$
\caln_\varepsilon=\{\rho\in W\mid \,  p_1(K_{\rho^0})-p_1(K_\rho)<\varepsilon \hbox{ and }
|L_\rho|- |L_{\rho^0}|<\varepsilon\} \ .
$$
Consider
the unique rotation $R_\rho\in SO(d)$ such that
$R_\rho(L_\rho)=-|L_\rho|e_1$ and
$R_\rho(q)=q$ if $q\in{\rm EV}(e_1,K_\rho)^\perp$
(if $\varepsilon$ is small enough, $L_\rho$ is not a positive
multiple of $e_1$ when $\rho\in\caln_\varepsilon$, thus $R_\rho$ is well defined).
If $\varepsilon$ is small enough, we check, as in \cite[Proof of Theorem~3.2]{Ha}
that the smooth maps
$\phi_-:\caln_\varepsilon\to(\bbr^{d-1})^k$ and $\phi_+:\caln_\varepsilon\to(\bbr^{d-1})^{m-k-2}$
given by
$$
\phi_-(\rho)=(P(\rho_1),\dots,P(\rho_k))
\ \hbox{ and } \
\phi_+(\rho)=\big( P(R_\rho(-\rho_{k+1})),\dots,P(R_\rho(-\rho_{m-1})) \big)
$$
are $O(d-1)$-equivariant and give rise to a $O(d-1)$-equivariant chart
$$
\phi=(\phi_-,\phi_+) : \caln_\varepsilon\to (\bbr^{d-1})^k\times (\bbr^{d-1})^{m-2-k} =
(\bbr^{d-1})^{m-2} \, ,
$$
where $(\bbr^{d-1})^n$ is endowed with the diagonal action of $O(d-1)$.
One has $\phi(\rho^0)=0$.
The subspaces
$$
D_+=\{\rho\in \caln_\varepsilon \mid \, |L_\rho|=|L_{\rho^0}| \}
\ \hbox{ and } \
D_-=\{\rho\in \caln_\varepsilon \mid \, K_\rho=K_{\rho^0}\}
$$
are submanifolds of dimensions $k(d-1)$ and $(m-k-2)(d-1)$ respectively,
satisfying $\phi(D_+)\subset (\bbr^{d-1})^k\times 0$ and
$\phi(D_+)\subset 0\times (\bbr^{d-1})^{m-2-k}$.

As in \cite[Proof of Theorem~3.2]{Ha}, we prove that $f$
restricts to Morse functions on $D_\pm$. The single critical point $\rho^0$
is a minimum on $D_+$ and a maximum on $D_-$. Therefore,
the Hessian form of $f\pcirc\phi^\mun$ is
positive definite on $T_0 (\bbr^{d-1})^k$ and
negative definite on $T_0 (\bbr^{d-1})^{m-2-k}$.
We have seen above that the $O(d-1)$-action on these subspaces is the standard diagonal
action. By \cite[Lemma~4.5]{Wn}, this implies the last assertion of \proref{Psurgery1}.
\cqfd

\begin{Proposition}\label{Psurgery2}
Suppose, in \proref{Psurgery1}, that $|J|=2$. Then
\begin{enumerate}
\item $\cha{m}{d}(\beta)=\cha{m}{d}(\alpha)\,\sharp\, (S^{d-1}\times S^{(m-3)(d-1)-1})$
\item $\nua{m}{3}(\beta)=\nua{m}{3}(\alpha)\,\sharp\, \overline{\bbc P}^{m-3}$
\end{enumerate}
\end{Proposition}

\preu
Since $\cals_m(\beta)=\cals_m(\alpha)\cup \{J\}$, the chamber
$\alpha$ is not  $\llangle{\{m,m-3,m-2,\dots,1\}}$ by \remref{111-rem}.
By \proref{connectivity},
$\cha{m}{d}(\alpha)$ is  ($d-2$)-connected. Hence, the sphere
$S^{d-2}\subset\cha{m}{d}(\alpha)$ on which the surgery of
\proref{Psurgery1} is performed is null-homotopic. We may assume that
$m\geq 4$ and $d\geq 2$ since \proref{Psurgery2} is empty for $m=3$
and $\cha{m}{1}(\alpha)=\cha{m}{1}(\beta)=\emptyset$ because of genericity.
Therefore $2(d-2)<{\rm dim\,}\cha{m}{d}(\alpha)$, from which we deduce that
$S^{d-2}\subset\cha{m}{d}(\alpha)$ is isotopic to a sphere
contained in a disk. Observe that we are dealing with stably parallelizable
manifolds (for instance, $\cha{m}{d}(-)$ is the pre-image of a regular value
of a map from a product of spheres to $\bbr^d$).
Part 1 then follows from standard results in surgery, see
e.g.~\cite[Proposition~11.2 and p.~188]{Ko}.

As for Part 2, we have
$$
\cha{m}{3}(\beta)\approx\big(\cha{m}{3}(\alpha)\psetminus\,
(S^1\times D^{2(m-3)})\big)
\cup_{S^1\times S^{2(m-3)-1}}
\big(D^2\times S^{2(m-3)-1}\big) \, .
$$
The quotient space $S^1\times D^{2(m-3)}$ by the action of $SO(2)$
is a disk $D^{2(m-3)}$. On the other hand, consider the tautological
line bundle $E\to\bbc P^{m-4}$, where
$E=\{(v,\ell)\in \bbc^{m-3}\times\bbc P^{m-4}\mid v\in\ell\}$.
Seeing $D^2$ as the unit disk in $\bbc$,
the map $g:D^2\times S^{2(m-3)-1}\to E$ given by
$g(z,w)=(zw,\bbc w)$ descends to an embedding from
$SO(2)\big\backslash (D^2\times S^{2(m-3)-1})$ to a neighborhood
of the zero section of $E$. It follows that $\cha{m}{3}(\beta)$
is diffeomorphic to $\cha{m}{3}(\alpha)$ blown up at one point,
which implies Assertion 2 (see e.g. \cite[pp.~214--216]{MDS}).
\cqfd

\begin{Example}\label{mdot2}
The chamber $\llangle{\{m,m-3,m-2,\dots,2\}}$. \ \rm
Let $\alpha=\llangle{\{m,m-3,m-2,\dots,2\}}$ and
$\beta=\llangle{\{m,m-3,m-2,\dots,1\}}$. Then
$\cals_m(\beta)=\cals_m(\alpha)\cup\{m,m-3,m-2,\dots,1\}$.
By \proref{Psurgery1}, $\cha{m}{d}(\beta)$
is obtained from $\cha{m}{d}(\alpha)$ by an $O(d-1)$-equivariant
surgery of index $(d-1)(m-3)-1$. Then, conversely,
$\cha{m}{d}(\alpha)$
is obtained from $\cha{m}{d}(\beta)$ by an $O(d-1)$-equivariant
surgery of index $d-2$. By \exref{111} 
and \proref{connectivity}, $\cha{m}{d}(\beta)
\approx (S^{d-1})^{m-3}\times S^{d-2}$
while $\cha{m}{d}(\alpha)$ is $(d-2)$-connected. This implies that the surgery on $\cha{m}{d}(\beta)$ is performed on
a tubular neighborhood of $pt\times S^{d-2}$.
Thus, Part 1 of \proref{Psurgery2} is not true, 
but one has
\begin{equation}\label{mdot2-eq1}
\cha{m}{d}(\alpha)
\approx
\big[\big((S^{d-1})^{m-3}\,\psetminus\,B\big)\times S^{d-2}
\big] \cup_{\partial B\times S^{d-2}}
(\partial B\times D^{d-1}) \, ,
\end{equation}
where $B$ is a $((m-3)(d-1))$-disk in $(S^{d-1})^{m-3}$. 
This is not a very simple expression,
except when $d=2$ where $\nua{m}{2}(\alpha)=\nua{m}{2}(\alpha)$ becomes
\begin{equation}\label{mdot2-eq2}
\nua{m}{2}(\alpha)
\approx
(S^{d-1})^{m-3}\,\sharp\,(S^{d-1})^{m-3} \, .
\end{equation}
On the other hand, Part 2 of \proref{Psurgery2} is valid and we get
\begin{equation}\label{mdot2-eq3}
\nua{m}{3}(\alpha)
\approx
(S^2)^{m-3}\,\sharp\,\overline{\bbc P}^{m-3} \, .
\end{equation}
It was observed by D.~Sch\"utz that there are only three chambers $\alpha\in\chrmcr$
such that $\cals_m(\alpha)$ contains $A=\{m,m-3,m-2,\dots,2\}$. These are 
\begin{equation}
\begin{array}{ccl}
\alpha   &=&  \llangle{A}  \\
\alpha'  &=&  \llangle{A,\{m,m-2\}}  \\
\alpha'' &=&  \llangle{A,\{m,m-1\}} \ .
\end{array}
\end{equation}
Indeed, if $A$ is short, then $\{m-1,m-2,1\}$ is long and one cannot add to $A$ 
a gene containing $3$ elements.
As $\cals_m(\alpha')=\cals_m(\alpha)\cup\{m,m-2\}$ and 
$\cals_m(\alpha'')=\cals_m(\alpha')\cup\{m,m-1\}$, the chain and polygon spaces for
$\alpha'$ and $\alpha''$ may  be obtain from the above using \proref{Psurgery2}. 

\end{Example}

\begin{Example}The chamber $\llangle{\{m,p\}}$. \ \rm
For $p\geq 2$, one has 
$\cals_m(\llangle{\{m,p\}})=\cals_m(\llangle{\{m,p-1\}})\cup \{m,p\}$ and 
$\cals_m(\llangle{\{m,1\}})=\cals_m(\llangle{m})\cup \{m,1\}$.
Using \proref{Psurgery2} and \exref{mmm}, one sees that 
\begin{center}
\begin{tabular}{c|c|c}
\multicolumn{3}{c}{\bf The chamber $\alpha=\llangle{\{m,p\}}$ }
\\ \hline\rule{0pt}{3ex}
 $\nua{m}{2}(\alpha)$ & $\nua{m}{3}(\alpha)$ & $\cha{m}{d}(\alpha)$
\\[1pt]
\hline\rule{0pt}{3ex}   
$p\,(S^{1}\times S^{m-4})$  &  $\bbc P^{m-3}\,\sharp\, p\,\overline{\bbc P}^{m-3}$ &
$p\,(S^{d-1}\times S^{(m-3)(d-1)-1})$
\\ \hline
\end{tabular}
\end{center}
Here, $p$ times a manifold $V$ means the connected sum of $p$ copies of $V$
(hence, a sphere if $p=0$).
A representative of $\llangle{\{m,p\}}$ is given by 
$$(\underbrace{1,\dots,1}_p,\underbrace{2,\dots,2}_{m-p-1},2m-p-5) \, .$$
The tables of \cite{HRWeb}
show that, for $m\leq 9$ (See Section \ref{Stables} below for $m\leq 6$),
this representative is $a_{\rm min}(\llangle{\{m,p\}})$, except
for $p=0,1$. 
As $\llangle{\{m,1\}}=\llangle{m-1}^+$, \proref{Paddtiny-nua}
gives the diffeomorphism
$$
\bbc P^{m-3}\,\sharp\,\overline{\bbc P}^{m-3}\approx
\nua{m}{3}(\llangle{\{m,1\}})\approx \nua{m}{3}(\llangle{m}^+)
\approx
S^2\times_{S^1} S^{2(m-3)-1} \, .
$$
In the case $m=5$, we get the two topological descriptions of the 
Hirzebruch surface (see, e.g. \cite[Ex.~6.4]{MDS}). 
\end{Example}

\section{Tables for $m=4,5,6$}\label{Stables}

For any $m$, there is the ``trivial'' chamber $\llangle{}$, where $a_m$ is so long that 
the corresponding chain or polygon spaces are empty. 
When $m=4$, Examples~\ref{mmm} and~\ref{111} give the remaining two chambers: 

\begin{center}
\begin{tabular}{rcc|c|c|c}
\multicolumn{5}{c}{\bf Table A: \ $m=4$}\\
\hline \rule{0pt}{3ex}
& $\alpha$ & $a_{\rm min}(\alpha)$ & 
\hspace{1pt} $\nua{4}{2}(\alpha)$ \kern 5mm &
$\nua{4}{3}(\alpha)$ & $\cha{4}{d}(\alpha)$
\\[3pt]\hline
\rule{0pt}{3ex}
1&$\langle\rangle$ & $(0,0,0,1)$ 
\ccvi{}{$\emptyset$}{$\emptyset$}{$\emptyset$}
\hline \rule{0pt}{3ex}
2&
$\langle 4\rangle$ & $(1,1,1,2)$  
\ccvi{}{$S^{2(d-1)-1}$}{$\bbc P^1$}{$S^1$}
\hline \rule{0pt}{3ex}
3&
$\langle 41\rangle$ & $(0,1,1,1)$  
\ccvi{}{$(S^{d-1})\times S^{d-2}$}{$S^2$}{$S^1\,\dot\cup\,S^1$}
\hline
\end{tabular}
\end{center}

Recall that the column $\cha{4}{d}(\alpha)$ does not contain all the $4$-chains,
only those for which $a_4\geq a_i$ for $I=1,2,3$. For example,
$\cha{4}{d}(1,1,1\varepsilon)\approx T^1 S^{d-1}$, the unit tangent bundle
to $S^{d-1}$. For a complete classification of $4$-chains, see~\cite{Ha}.

\sk{2}
When $m=5$, there are seven chambers. Lines 2 and 7 come from Examples~\ref{mmm} and~\ref{111}.
The symbols $\Sigma_g$ denotes the orientable surface of genus $g$
and $T^r$ is the torus $(S^1)^r$.
Within the central block, each line is obtained from the previous one by \proref{Psurgery2}.

\begin{center}
\begin{tabular}{rcc|c|c|c}
\multicolumn{5}{c}{\bf Table B: \ $m=5$}\\
\hline \rule{0pt}{3ex}
& $\alpha$ & $a_{\rm min}(\alpha)$ & 
\hspace{1pt} $\nua{5}{2}(\alpha)$ \kern 5mm &
$\nua{5}{3}(\alpha)$ & $\cha{5}{d}(\alpha)$
\\[3pt]\hline
\rule{0pt}{3ex}
1&$\langle\rangle$ & $(0,0,0,0,1)$ 
\ccvi{}{$\emptyset$}{$\emptyset$}{$\emptyset$}
\hline \rule{0pt}{3ex}
2&
$\langle 5\rangle$ & $(1,1,1,1,3)$  
\ccvi{}{$S^{3(d-1)-1}$}{$\bbc P^2$}{$S^2$}
3&
$\langle 51\rangle$ & $(0,1,1,1,2)$
\ccvi{}{$S^{d-1}\times S^{2(d-1)-1}$}{$\bbc P^2\,\sharp\,\overline{\bbc P}^2$}{$T^2$}
4&
$\langle 52\rangle$ & $(1,1,2,2,3)$
\ccvi{}{$2\,[S^{d-1}\times S^{2(d-1)-1}]$}{$\bbc P^2\,\sharp\,2\,\overline{\bbc P}^2$}{$\Sigma_2$}
5&
$\langle 53\rangle$ & $(1,1,1,2,2)$ 
\ccvi{}{$3\,[S^{d-1}\times S^{2(d-1)-1}]$}{$\bbc P^2\,\sharp\,3\,\overline{\bbc P}^2$}{$\Sigma_3$}
6&
$\langle 54\rangle$ & $(1,1,1,1,1)$
\ccvi{}{$4\,[S^{d-1}\times S^{2(d-1)-1}]$}{$\bbc P^2\,\sharp\,4\,\overline{\bbc P}^2$}{$\Sigma_4$}
\hline \rule{0pt}{3ex}
7&
$\langle 521\rangle$ & $(0,0,1,1,1)$  
\ccvi{}{$(S^{d-1})^2\times S^{d-2}$}{$S^2\times S^2$}{$T^2\,\dot\cup\,T^2$}
\hline
\end{tabular}
\end{center}

\rm
Line 3 in Table B together with 
Equation~\equref{mdot2-eq3}, re-proves the classical
fact that 
$(S^2\times S^2)\,\sharp\,\overline{\bbc P}^2$ is diffeomorphic to 
$\bbc P^2\sharp\, 2\,\overline{\bbc P}^2$.

\sk{2}
When $m=6$, there are 21 chambers. In order to save space, we did not give $a_{\rm min}(\alpha)$
(they can be found in \cite[Table~6]{HR}). The first line of each block is obtained from
\S~\ref{atiny}--\ref{SSurg}. Then,
each line is obtained from the previous one by \proref{Psurgery2}.

\renewcommand{\pvi}[7]{\stepcounter{compol}{\small\thecompol} & %
$\langle #7\rangle$  }

\newpage
{\footnotesize
\hskip -10mm
\setcounter{compol}{0}
\begin{supertabular}{rl|c|c|c}
\multicolumn{5}{c}{\bf Table C: \ $m=6$}\\
\hline \rule{0pt}{3ex}
& $\alpha$ & 
\hspace{1pt} $\nua{6}{2}(\alpha)$ \kern 5mm &
$\nua{6}{3}(\alpha)$ & $\cha{6}{d}(\alpha)$
\\[3pt]\hline
\rule{0pt}{3ex}
\pvi{0}{0}{0}{0,0,0,0,0,1}{1}{1}{}%
\ccvi{}{$\emptyset$}{$\emptyset$}{$\emptyset$}%
\hline \rule{0pt}{3ex}
\pvi{1}{1}{1}{1,1,1,1,1,4}{9}{2}{6}%
\ccvi{}{$S^{4(d-1)-1}$}{$\bbc P^3$}{$S^3$}%
\pvi{2}{2}{2}{0,1,1,1,1,3}{7}{3}{61}%
\ccvi{$S^{d-1}\!\times\! S^{3(d-1)-1}$}{$R(d)$ \fno{1}}%
{$\bbc P^3\,\sharp\,\overline{\bbc P}^3$}{$S^1\!\times\! S^2$}%
\pvi{3}{3}{3}{1,1,2,2,2,5}{13}{4}{62}
\ccvi{}{$2R(d)$}%
{$\bbc P^3\,\sharp\,2\,\overline{\bbc P}^3$}{$2(S^1\!\times\! S^2)$}%
\pvi{4}{4}{4}{1,1,1,2,2,4}{11}{6}{63}
\ccvi{}{$3R(d)$}%
{$\bbc P^3\,\sharp\,3\,\overline{\bbc P}^3$}{$3(S^1\!\times\! S^2)$}%
\pvi{5}{5}{5}{1,1,1,1,2,3}{9}{10}{64}
\ccvi{}{$4R(d)$}%
{$\bbc P^3\,\sharp\,4\,\overline{\bbc P}^3$}{$4(S^1\!\times\! S^2)$}%
\pvi{6}{6}{6}{1,1,1,1,1,2}{7}{12}{65}
\ccvi{}{$5R(d)$}%
{$\bbc P^3\,\sharp\,5\,\overline{\bbc P}^3$}{$5(S^1\!\times\! S^2)$}%
\hline \rule{0pt}{3ex}
\pvi{3}{2}{0}{0,0,1,1,1,2}{5}{5}{621}
\ccvi{}{$S^{d-1}\!\times\! \cha{5}{d}(\llangle{51})$ \fno{2}}%
{$S^2\!\times\!_{S^1}(S^2\!\times\! S^3)$}{$T^3$}%
\pvi{4}{3}{1}{1,1,2,3,3,5}{15}{14}{621,63}
\ccvi{}{$[S^{d-1}\!\times\! \cha{5}{d}(\llangle{51})]\,\sharp\,R(d)$}%
{$[S^2\!\times\!_{S^1}(S^2\!\times\! S^3)]\,\sharp\,\overline{\bbc P}^3$}{$T^3\sharp\,(S^1\!\times\!S^2)$}%
\pvi{5}{4}{2}{1,1,2,2,3,4}{13}{15}{621,64}
\ccvi{}{$[S^{d-1}\!\times\! \cha{5}{d}(\llangle{51})]\,\sharp\,2R(d)$}{$[S^2\!\times\!_{S^1}(S^2\!\times\! S^3)]\,\sharp\,2\,\overline{\bbc P}^3$}{$T^3\sharp\,2(S^1\!\times\! S^2)$}%
\pvi{6}{5}{3}{1,1,2,2,2,3}{11}{16}{621,65}
\ccvi{}{$[S^{d-1}\!\times\! \cha{5}{d}(\llangle{51})]\,\sharp\,3R(d)$}{$[S^2\!\times\!_{S^1}(S^2\!\times\! S^3)]\,\sharp\,3\,\overline{\bbc P}^3$}%
{$T^3\sharp\,3(S^1\!\times\! S^2)$}%

\hline \rule{0pt}{3ex}
\pvi{4}{2}{0}{0,1,1,2,2,3}{9}{7}{631}
\ccvi{}{$S^{d-1}\!\times\! \cha{5}{d}(\llangle{52})$  \fno{2}}%
{$S^2\!\times\!_{S^1}2(S^2\!\times\! S^3)$}{$\Sigma_2\!\times\! S^1$}%
\pvi{5}{3}{1}{1,2,2,3,4,5}{17}{17}{631,64}
\ccvi{}{$[S^{d-1}\!\times\! \cha{5}{d}(\llangle{52})]\,\sharp\,R(d)$}%
{$S^2\!\times\!_{S^1}2(S^2\!\times\! S^3)\sharp\,\overline{\bbc P}^3$}%
{$(\Sigma_2\!\times\! S^1)\,\sharp\,(S^1\!\times\! S^2)$}%
\pvi{6}{4}{2}{1,2,2,3,3,4}{15}{18}{631,65}
\ccvi{}{$[S^{d-1}\!\times\! \cha{5}{d}(\llangle{52})]\,\sharp\,2R(d)$}%
{$S^2\!\times\!_{S^1}2(S^2\!\times\! S^3)\sharp\,2\,\overline{\bbc P}^3$}%
{$(\Sigma_2\!\times\! S^1)\,\sharp\,2(S^1\!\times\! S^2)$}%

\hline \rule{0pt}{3ex}
\pvi{5}{2}{0}{0,1,1,1,2,2}{7}{11}{641}
\ccvi{}{$S^{d-1}\!\times\! \cha{5}{d}(\llangle{53})$  \fno{2}}%
{$S^2\!\times\!_{S^1}3(S^2\!\times\! S^3)$}{$\Sigma_3\!\times\! S^1$}%
\pvi{6}{3}{1}{1,2,2,2,3,3}{13}{21}{641,65}
\ccvi{}{$[S^{d-1}\!\times\! \cha{5}{d}(\llangle{53})]\,\sharp\,R(d)$}%
{$S^2\!\times\!_{S^1}3(S^2\!\times\! S^3)\sharp\,\overline{\bbc P}^3$}%
{{$(\Sigma_3\!\times\! S^1)\,\sharp\,(S^1\!\times\! S^2)$}}%

\hline \rule{0pt}{3ex}
\pvi{6}{2}{0}{0,1,1,1,1,1}{5}{13}{651}
\ccvi{}{$S^{d-1}\!\times\! \cha{5}{d}(\llangle{54})$  \fno{2}}%
{$S^2\!\times\!_{S^1}4(S^2\!\times\! S^3)$}{$\Sigma_4\!\times\! S^1$}%

\hline \rule{0pt}{3ex}
\pvi{3}{0}{0}{0,0,0,1,1,1}{3}{9}{6321}%
\ccvi{}{$(S^{d-1})^{3}\!\times\! S^{d-2}$}{$(S^2)^3$}{$T^3\,\dot\cup\,T^3$}
\pvi{4}{1}{1}{1,1,1,3,3,4}{13}{8}{632}
\ccvi{}{$\cha{6}{d}(\llangle{632})$  \fno{3}}%
{$(S^2)^3\,\sharp\,\overline{\bbc P}^3$}{$2\,T^3$}%
\pvi{5}{2}{2}{1,1,1,2,3,3}{11}{19}{632,64}
\ccvi{}{$\cha{6}{d}(\llangle{632})\,\sharp\, R(d)$}%
{$(S^2)^3\sharp\,2\,\overline{\bbc P}^3$}%
{$2T^3\sharp\,(S^1\!\times\! S^2)$}%
\pvi{6}{3}{3}{1,1,1,2,2,2}{9}{20}{632,65}
\ccvi{}{$\cha{6}{d}(\llangle{632})\,\sharp\, 2R(d)$}%
{$(S^2)^3\sharp\,3\,\overline{\bbc P}^3$}%
{$2T^3\sharp\,2(S^1\!\times\! S^2)$}%
\hline\\
\multicolumn{5}{l}{\fno{1} $R(d)=S^{d-1}\!\times\! S^{3(d-1)-1}$. \hskip 12mm \fno{2}  see Table~B.  
\hskip 14mm \fno{3} See \exref{mdot2} }\\[2mm]
\end{supertabular}
}
\rm
The list for $\nua{6}{2}(\alpha)$ is present in \cite{Wa} with some short-hand justification.
A version of the column for $\nua{6}{3}(\alpha)$ is in \cite{HR}.

\sk{1}
For $m\geq 7$, the above procedure fails to give all the chambers, since
surgeries of higher index are needed. For example, for $m=7$, only
49 chambers out of 135 are reached.


\end{document}